\documentclass{amsart}

\usepackage{amsmath}
\usepackage{amsfonts}
\usepackage{amssymb}
\usepackage{amsthm}
\usepackage{comment}
\usepackage{epsfig}
\usepackage{psfrag}
\usepackage{mathrsfs}
\usepackage{amscd}
\usepackage[all]{xy}
\usepackage{rotating}
\usepackage{lscape}
\usepackage{amsbsy}
\usepackage{verbatim}
\usepackage{moreverb}

\newtheorem{theorem}{Theorem}[section]

\theoremstyle{definition}

\theoremstyle{remark}
\newtheorem{remark}[theorem]{Remark}
\newtheorem{example}[theorem]{Example}

\DeclareMathOperator{\holim}{holim}

\DeclareMathOperator{\Map}{Map}

\newcommand{\tmf}{\mathit{tmf}}

\newcommand{\CP}{\mathbb{CP}^\infty}

\newcommand{\GL}{\mathit{GL}}
\newcommand{\gl}{\mathit{gl}}

\def\cL{\mathcal L}
\def\cO{\mathcal O}

\def\CC{\mathbb C}
\def\GG{\mathbb G}

\def\SS{\mathbb S}

\def\ZZ{\mathbb Z}

\begin{document}

\title{
Lecture by Michael Hopkins: the string orientation of $\tmf$}
\author{{\rm Held on March $29^\text{th}$ 2007, at the Talbot workshop}\\\\Typed by A. Henriques}

\maketitle

There are a lot of aspects to the question of orientations.
There are things like why you would even look for such a thing, and what it means to have one.
Once you want one, there are the tricks of the trade used to produce it.
The latter involve some computations that might be yet another aspect...

\[
\sim\qquad\sim\qquad\sim
\]\phantom{.}

So our goal is to construct an $E_\infty$ map
\begin{equation}\label{MOtmf}
MO\langle8\rangle\to \tmf.
\end{equation}
Here, $MO\langle8\rangle$ is a synonym of $M\mathit{String}$, and stands for the Thom spectrum of the 7-connected cover of $BO$.
Recall that the connected covers of $BO$ fit in a tower
\[
\ZZ\times BO
\leftarrow BO
\leftarrow B\mathit{Spin}
\leftarrow B\mathit{String}.
\]
The space $BO$ is obtained from $\ZZ\times BO$ by killing its $\pi_0$.
$B\mathit{Spin}$ is obtained from $BO$ by killing $\pi_1(BO)=\ZZ/2$.
Finally, $B\mathit{String}$ is obtained from $B\mathit{Spin}$ by killing its first non-vanishing
homotopy group, namely
\[
\pi_4(B\mathit{Spin})=\ZZ.
\]
The first Pontryagin class of the universal vector bundle over $B\mathit{Spin}$ being twice the generator of
$H^4(B\mathit{Spin})= \ZZ$, that generator is usually called $\frac{p_1}{2}$.
The space $B\mathit{String}$ is the fiber of the map $\frac{p_1}{2}:B\mathit{Spin}\to K(\ZZ,4)$.
It was previously called $BO\langle8\rangle$;
its current optimistic name is due to connections with string theory.

The group $\pi_{2n}(MO\langle8\rangle)$ is the bordism group of string manifolds of dimension $2n$. 
Recall that $\pi_{2n}(\tmf)$ has a natural maps to $MF_n$, the group of modular forms of weight $n$.
At the level of homotopy groups, the map (\ref{MOtmf}) should then send a string manifold $M$ to its Witten genus
$\phi_W(M)$
\[
\begin{split}
\pi_{2n}(MO\langle8\rangle)&\rightarrow \pi_{2n}(\tmf) \rightarrow \mathit{MF}_n\\
[M]\hspace{.5cm}&\hspace{.8cm}\mapsto \hspace{1cm}\phi_W(M).
\end{split}
\]

The cohomology theory $\tmf$ was known to us [actually its 2-completion, and under the name $eo_2$] 
before we new about the connection with modular forms.
We were then looking for the map (\ref{MOtmf}) because of the following cohomology calculation.
The cohomology of $\tmf$ at the prime $2$ is the cyclic module 
$A/\langle\mathit{Sq}^1,\mathit{Sq}^2,\mathit{Sq}^4\rangle$ over the Steenrod algebra $A$.
That same cyclic module occurs as a sumand in the cohomology of $MO\langle8\rangle$, which
lead us to believe that there should be such a map.
The connection with elliptic curves was made while trying to make the map $MO\langle8\rangle\to\tmf$.
Indeed, producing a map from $MO\langle8\rangle$ into a complex oriented cohomology theory $E$
is something that one can do easily if the formal group associated to $E$ comes from an elliptic curve.
So the whole story of $\tmf$ had to do with that orientation.
It is only in retrospect that we noted that the map (\ref{MOtmf}) reproduces the Witten genus.

One interesting fact is that the map $\pi_{2n}(\tmf) \rightarrow MF_n$ is not quite an isomorphism.
It is an isomorphism after tensoring with $\ZZ[\frac{1}{2},\frac{1}{3}]$,
but it contains some torsion in its kernel, and its image is only a subgroup of finite index.
So that way, one learns some things about the Witten genus that one might not have known before.
For example, looking at $\pi_{24}$, we see a map
\begin{equation}\label{MF12}
\pi_{24}(\tmf)\simeq \ZZ\oplus \ZZ\longrightarrow \mathit{MF}_{12}\simeq \ZZ\oplus \ZZ.
\end{equation}
The group $MF_{12}$ has two generators $c_4^3$ and $\Delta$,
and it is interesting to note that the image of (\ref{MF12}) 
is the subgroup generated by $c_4^3$ and $24\Delta$.
So this gives a restriction on the possible values of the Witten genus.
Translating it back into geometry,
it says for example that on a $24$ dimensional string manifold, the index of the Dirac operator
with coefficients in the complexified tangent bundle is always divisible by $24$
\begin{equation}\label{div24}
\widehat{A}\big(M;T_\CC\big)\equiv 0 \qquad (\text{mod}\hspace{.2cm}24).
\end{equation}
The torsion in $\pi_*(\tmf)$ also give interesting ``mod 24 Witten genera", 
which are analogs of the mod $2$ indices of the Dirac operator in $KO$-theory 
(those exist for manifolds of dimension $1$ and $2$ mod $8$).
These are all facts for which there is no explanation in terms of index theory, or even string theory.

In short, one reason for wanting (\ref{MOtmf}) is that one gets more refined geometric information
about the Witten genus of string manifolds.
According to string theory, the Witten genus is supposed to be the index of the Dirac operator on the loop space $LM$ of $M$.
There ought to be some kind of structure on $LM$ that accounts for the factor of $24$ in (\ref{div24}),
but up to now, there is no explanation using the geometric approach.

\[
\sim\qquad\sim\qquad\sim
\]\phantom{.}

We now explain why $E_\infty$ maps out of $MO\langle8\rangle$,
have anything to do with elliptic curves.
In some sense, there is a very natural reason to expect a map like (\ref{MOtmf}).
To simplify the analysis, we work with the complex analog of $MO\langle8\rangle$, namely $MU\langle6\rangle$.

Considering the various connected covers of $BU$, one gets a tower of spaces
\[
\ZZ\times BU
\leftarrow BU
\leftarrow BSU
\leftarrow BU\langle6\rangle
\]
whose last term is the fiber of $c_2:BSU \to K(\ZZ,4)$.
These have companion bordism theories,
and $MU\langle6\rangle$ is the one corresponding to $BU\langle6\rangle$.
In other words, $MU\langle6\rangle$
is the Thom spectrum of the universal bundle over $BU\langle6\rangle$.

We recall Quillen's theorem in its formulation via formal groups.
Roughly speaking, it says that multiplicative maps from $MU$ into an (even periodic) complex orientable cohomology theory $E$
correspond to coordinates on the formal group $G:=\mathrm{spf}(E^0(\CP))$ associated to $E$.
In fact, the above statement is not quite accurate.
It is true that a coordinate gives you a map $MU\to E$, but the latter encode slightly more data.

To understand the subtlety, we begin with an analogy.
Multiplicative maps from the suspension spectrum of $BU$ into $E$ also correspond to some structures on $G$.
The important thing about $BU$ is that there is an inclusion $\CP\hookrightarrow BU$ 
exhibiting $E_*(BU)$ as the free commutative $E_*$-algebra on the $E_*$-module $E_*(\CP)$.
The sloppy analysis goes as follows.
A multiplicative map
\begin{equation}\label{BUtoE}
\Sigma^\infty BU_+ \to E
\end{equation}
corresponds to a ring homomorphism
\[
E_*(BU)\to E_*,
\]
and since $E_*(BU)$ is the symmetric algebra on $E_*(\CP)$, those correspond to $E_*$-module maps
\[
E_*(\CP)\to E_*.
\]
The latter are then elements of $E^0(\CP)$, in other words functions on $G$.
The problem in the above reasoning is that
if we want (\ref{BUtoE}) to be a multiplicative map,
the base point has to go to $1$.
So the base point of $\CP$ has to go to $1\in E$,
and therefore we don't get all functions $f$ on $G$, but only those satisfying $f(e)=1$, where 
$e\in G$ is the unit.
Now if we run all this through the Thom isomorphism, we find that multiplicative maps
\[
MU\to E
\]
are expressions of the form $f(z)/dz$ on $G$, such that the residue of $dz/f(z)$ is $1$.

Let $\cO$ denote the structure sheaf of $G$, let $e:S\to G$ be the identity section, 
and let $p:G\to S$ be the projection of $G$ onto the base scheme $S:=\mathrm{spec}(\pi_0 E)$.
Let also $\cO(-e)$ be the line bundle over $G$, whose sections are functions vanishing at zero.
Expressions of the form $f(z)/dz$ can then be understood as sections of the line bundle
\begin{equation}\label{OpeO}
\cO(-e)\otimes p^* e^*\cO(-e)^{-1}.
\end{equation}
The fiber of (\ref{OpeO}) over the origin is canonically trivialized, and
the condition $\mathrm{res}(dz/f(z))=1$ means that the value at $e$ of our section 
is equal to $1\in e^*\big(\cO(-e)\otimes p^* e^*\cO(-e)^{-1}\big)$.

Now, we would like to do something analogous for $BSU$ instead of $BU$.
Let $L$ denote the tautological line bundle over $\CP$.
The thing that allowed us to do the previous computation was the map
\[
\CP\to BU
\]
classifying $1-L$.
That map doesn't lift to $BSU$ because $c_1(1-L)\not =0$,
but the problem goes away as soon as one takes the product of two such bundles.

Let $L_1$ and $L_2$ denote the line bundles over $\CP\times\CP$ given by $L_1:=L\times 1$, and $L_2:=1\times L$ respectively.
Since $c_1((1-L_1)\otimes(1-L_2))=0$, we get a map
\[
\CP\times\CP\xrightarrow{(1-L_1)\otimes(1-L_2)}BSU.
\]
If $E$ is a complex orientable cohomology theory,
then we find that multiplicative maps
\[
\Sigma^\infty BSU_+ \to E,
\]
equivalently ring homomorphisms
\[
E_*(BSU)\to E_*,
\]
give rise to functions $f(x,y)$ on $G\times G$ satisfying
\begin{align}
\notag f(e,e)&=1,\\
\label{sym2co} f(x,y)&=f(y,x),\\
\notag \text{and}\qquad\quad f(y,z)f(x,y+_Gz)&=f(x,y)f(x+_Gy,z).\hspace{1cm}
\end{align}
In other words, they are functions on
$G$ with values in the multiplicative group,
that are rigid, and that are symmetric 2-cocycles. 
The last condition in (\ref{sym2co}) is obtained by expanding the virtual bundle
$(1-L_1)\otimes(1-L_2)\otimes(1-L_3)$ over $(\CP)^3$ in the following two ways:
\[
\begin{split}
(1-L_1)(1-L_2)(1-L_3)&=(1-L_1)(1-L_3)+(1-L_1)(1-L_2)-(1-L_1)(1-L_2L_3)\\
&=(1-L_1)(1-L_3)+(1-L_2)(1-L_3)-(1-L_1L_2)(1-L_3),
\end{split}
\]
from which we deduce that
\[
(1-L_2)(1-L_3) + (1-L_1)(1-L_2L_3) = (1-L_1)(1-L_2) + (1-L_1L_2)(1-L_3).
\]
In fact, the conditions (\ref{sym2co}) characterize
homomorphisms from $E_*(BSU)$ into any $E_*$-algebra.

In the case of $BU\langle6\rangle$,
multiplicative maps
\[
\Sigma^\infty BU\langle6\rangle_+ \to E,
\]
equivalently ring homomorphisms
\[
E_*(BU\langle6\rangle)\to E_*,
\]
give rise to functions of three variables $f:G^3\to\GG_m$ that satisfy the following conditions:
they are rigid, meaning that $f(0,0,0)=1$,
they are symmetric,
and they are 2-cocycles in any two of the three variables.

\begin{remark}
That kind of analysis stops at $BU\langle6\rangle$ because it
is the last connected cover of $BU$ with only even dimensional cells.
\end{remark}

\begin{remark}
There is an interesting analogy with classical group theory.
Let $\Gamma$ be a group, and let $I\subset \ZZ[\Gamma]$ denote its augmentation ideal.
A function $I\to A$ is the same as a function $f:G\to A$ satisfying $f(e)=0$.
A function $I^2\to A$ is the same thing as a function on $G\times G$, that is rigid, and a symmetric 2-cocycle.
Finally, a function $I^3\to A$ is the same thing as a function $G^3\to A$
that is rigid, symmetric, and a 2-cocycle in any two of the three variables.
So, in some sense, the connected covers of $BU$ correspond to
taking powers of the augmentation ideal of the group ring of a formal group.
\end{remark}

We now proceed to $MU\langle6\rangle$.
Recall from our discussion with $MU$ that doing the Thom isomorphism amounts to replacing functions by sections of a line bundle.
Since multiplicative maps $BU\langle6\rangle\to E$ corresponded to functions on $G^3$,
we expect multiplicative maps $MU\langle6\rangle\to E$ to produce sections of some line bundle $\Theta$ over $G^3$.
Letting $\cL:=\cO(-e)$,
we can describe $\Theta$ by writing down its fiber $\Theta_{(x,y,z)}$ over a point\footnote{$G^3$ is a formal scheme, and so it doesn't have many points.
To make sense of (\ref{TH}) and (\ref{sym2coc}), one can use the formalism of ``functor of points'' and work with maps $T\to G^3$, that one views as $T$-parametrized families of points.
}
$(x,y,z)\in G^3$.
It is given by
\begin{equation}\label{TH}
\Theta_{(x,y,z)}:=\frac{\hspace{.3cm}\cL_{x+y+z}\otimes\cL_x\otimes\cL_y\otimes\cL_y\hspace{.3cm}}
{\cL_{x+y}\otimes\cL_{x+z}\otimes\cL_{y+z}\otimes\cL_e},
\end{equation}
where $+$ denotes to the operation of $G$.
Multiplicative maps
\[
MU\langle6\rangle\to E
\]
then correspond to sections $s$ of $\Theta$ that are rigid, symmetric, and $2$-cocycles in any two 
of the three variables.
These conditions make sense because the two sides of each one of the equations 
\begin{equation}\label{sym2coc}
\begin{split}
s(e,e,e)&=1\\
s(x,y,z)=s(&y,x,z)=s(x,z,y)\\
s(y,z,v)s(x,y+z,v)&=s(x,y,v)s(x+y,z,v),
\end{split}
\end{equation}
are sections of (canonically) isomorphic line bundles.
For example, the first of the above equations makes sense because
the fiber $\Theta_{(e,e,e)}$ is canonically trivialized.

Now here's the thing that was inspiring to us:
if $J$ is an elliptic curve, then the line bundle $\Theta$ is trivial over $J^3$.
This is a special case of a general fact for abelian varieties called the theorem of the cube.

To understand why $\Theta$ is a trivial line bundle,
recall that given a divisor $D$ on $J$ of degree zero, there is a meromorphic function $f$
with that given divisor iff the points of $D$ add up to zero (taking multiplicities into account).
In particular, given two points $x,y$ on our elliptic curve $J$, there exists a meromorphic function $f$
with simple poles at $-x$ and $-y$, and a simple zeroes at $-x-y$ and $e$.
But that function is only well defined up to scale, and there is no canonical choice for it.
In other words, the line bundle
\begin{equation}\label{O4}
\cO\big(\!-[-x-y]+[-x]+[-y]-[e]\,\big),
\end{equation}
is trivial, but not trivialized.
On the other hand, if we divide (\ref{O4}) by its fiber over zero,
then it acquires a canonical trivialization.
Fix points $x,y\in J$, and
consider the restriction of $\Theta$ to the subscheme $\{x\}\times\{y\} \times J$.
We then have a canonical isomorphism between $\Theta|_{\{x\}\times\{y\} \times J}$
and the quotient of (\ref{O4}) by its value at zero.
So we get canonical trivializations of each such restriction.
These trivializations then assemble to a trivialization of $\Theta$.

Note that $\Theta$ is not just trivial, it is canonically trivialized.
Therefore it makes sense to talk about the section ``1'' of $\Theta$.
If we take that section,
and pull it back along any map $J^m\to J^3$ then we'll always get the
section ``1'' of (another) canonically trivial line bundle. 
So any conditions one might decide to impose on a section of $\Theta$, 
for example the conditions (\ref{sym2coc}), will be automatically satisfied by our distinguished section ``1''.

The consequence of the above discussion is that if the formal group of $E$ comes from an elliptic curve,
then we get a canonical solution to the equations (\ref{sym2coc}). 
In particular, we get a canonical multiplicative map
\[
MU\langle6\rangle\to E.
\]
Now, there is an analog of that for $MO\langle8\rangle$ which involves adding one more condition
to the list (\ref{sym2coc}).
That condition is automatically satisfied for the same reasons as above,
and one finds that there is a canonical multiplicative map 
\begin{equation}\label{MO8toE}
MO\langle8\rangle\to E
\end{equation}
as soon as the formal group of $E$ comes from an elliptic curve.

If $J$ is an elliptic curve defined over a ring $R$, and $\varphi:R\to R'$ is a ring homomorphism,
one gets an induced elliptic curve $J'$ over $R'$.
Let $E$, $E'$ be complex orientable cohomology theories whose associated formal groups
are the formal completions of $J$ and $J'$.
If $f:E\to E'$ is a map of spectra with $\pi_0(\varphi)=f$,
then the maps (\ref{MO8toE}) induce a homotopy commutative diagram
\begin{equation}\label{MO8EandE'}
\begin{matrix}
\xymatrix@R=.5cm{
&MO\langle8\rangle\ar[dr]\ar[dl]&\\
E\ar[rr]|{\,f\,}&&E'.
}
\end{matrix}
\end{equation}
This is what led to the idea of assembling all cohomology theories coming from elliptic curves\footnote{
For the statement (\ref{tmf=}) to actually be correct, one needs to include the multiplicative group $\GG_m$, as well as the additive group $\GG_a$, in our definition of ``elliptic curves''.}
into a single cohomology theory
\begin{equation}\label{tmf=}
\tmf = \underset{\text{elliptic curves $J$}}{\holim}\Big( \text{cohomology theory $E$ associated to $J$}\Big).
\end{equation}
The maps (\ref{MO8EandE'}) then assemble into a map $MO\langle8\rangle\to \tmf$.
That map reproduces the Witten genus at the level of homotopy groups,
and is then an explanation of why the Witten genus of a string manifold is a modular form.

\[
\sim\qquad\sim\qquad\sim
\]\phantom{.}

So far, we have addressed the questions of why one should be interested in a map like (\ref{MOtmf}),
and why one could expect there to be one.
We now describe a homotopy theoretic way of producing that map.

To make the previous arguments actually work, one would need to do things in a much more rigid way.
Indeed, the maps (\ref{MO8toE}) are only multiplicative up to homotopy, and the triangles (\ref{MO8EandE'}) only commute up to homotopy.
To get a map of spectra $MO\langle8\rangle\to \tmf$, one would need to rigidify the triangles (\ref{MO8EandE'}).
And then, to get it to be a map of $E_\infty$ ring spectra, one would need to upgrade the maps 
(\ref{MO8toE}) to $E_\infty$ maps.
That shall not be our strategy.
Instead, we will produce the map (\ref{MOtmf}) directly, by using the decomposition of $\tmf$ into its various $p$-completions
and $K(n)$-localizations.
Of course, using derived algebraic geometry, there are now more conceptual ways of producing that map 
[J. Lurie interrupts and claims that derived algebraic geometry cannot produce the map].

Another example of something that one can construct using the same methods, is the map from spin bordism to $KO$ theory
\begin{equation}\label{MKO}
M\mathit{Spin}\to KO
\end{equation}
that sends a spin manifold to its $\widehat A$-genus, namely the index of the Dirac operator.
Equivalently, (\ref{MKO}) is the $KO$-theory orientation of spin bundles.
Note that before our techniques, there was no way known of producing the map (\ref{MKO}) using the methods of homotopy theory:
one had to use geometry.
The construction of (\ref{MKO}) will be our warm-up, before trying to get the more interesting one
\[
MO\langle8\rangle\to \tmf.
\]

Given a vector bundle $\zeta$ over a space $X$, we let $X^\zeta$ denote the associated Thom space.
One should think of $X^\zeta$ as a twisted form of $X$.
Actually, if we view $X^\zeta$ as a spectrum, then we should rather say that $X^\zeta$ is a twisted form of 
the suspension spectrum $\Sigma^\infty X_+$.
From now on, we shall abuse notation, and write $X\wedge-$ instead of $\Sigma^\infty X_+\wedge-$.

Similarly, given a cohomology theory $E$, the spectrum $X^\zeta\wedge E$ is a twisted form of $X\wedge E$.
An $E$-orientation is then a trivialization of the bundle of spectra obtained by fiberwise smashing $E$ with the (compactified) fibers of $\zeta$.
The Thom isomorphism is the induced equivalence between $X^\zeta\wedge E$ and $X\wedge E$.
In fact, we didn't need to start with a vector bundle: a spherical fibration is enough to produce a Thom spectrum.

If $E$ is not multiplicative, then there are no further conditions that one could impose.
But if $E$ is $E_\infty$, then we can ask that the above isomorphisms be isomorphisms of $E$-modules.
That's the concept of an $E_\infty$ orientation.
One reason for looking for $E_\infty$ orientation instead of just orientations, is that it simplifies the
computations.

If $R$ is a spectrum equipped with a homotopy associative product $R\wedge R\to R$,
we define $\GL_1(R)\subset \Omega R$ to be the subset of the zeroth space of $R$ where we only take
the connected components corresponding to the units of $\pi_0(R)$.
It satisfies
\[
[X,\GL_1(R)]=R^0(X)^\times
\]
for all unbased spaces $X$.
If $X$ has a base point, then the above equation is not quite correct.
Since the base point of $X$ has to go to the base point $1\in\GL_1(R)$,
the better way to write this is
\[
[X,\GL_1(R)]=\big(1+\widetilde{R}^0(X)\big)^\times.
\]
In words, it is the group of invertible elements of $R^0(X)$ that restrict to 1 at the base point.
If $R$ is an $A_\infty$-ring spectrum, then $\GL_1(R)$ is an $A_\infty$-space, and thus a loop space. 
We then let $B\GL_1(R)$ denote its classifying space.
\begin{example} $B\GL_1(\SS)$ is the classifying space for spherical fibrations.
\end{example}

If we have a map
\[
\zeta: X\to B\GL_1(\SS),
\]
then we get a spherical fibration over $X$, and thus a corresponding Thom spectrum $X^\zeta$.
Now, if we start instead with a map
\[
\zeta: X\to B\GL_1(R),
\]
then we can form an analogous construction in the world of $R$-modules.
Let $P$ be the $\GL_1(R)$-principal bundle associated to $\zeta$
\[
\xymatrix{
\GL_1(R)\ar[r]&P\ar[d]&\\
&X\ar[r]^(.35)\zeta&B\GL_1(R).
}
\]
We then define
\[
X^\zeta:=\Sigma^\infty P_+\underset{\Sigma^\infty \GL_1(R)_+}{\wedge}R.
\]
That construction can be understood as follows.
Each one of the fibers of $P$ is a copy of $\GL_1(R)$, more precisely a torsor for $\GL_1(R)$.
The operation $-\wedge_{\Sigma^\infty \GL_1(R)_+}\!R$ then converts that torsor into a copy of $R$.
So for each point of $X$, one gets a copy of $R$.
If $\zeta$ is the trivial map, one has $X^\zeta=X\wedge R$.
So for general $\zeta$, the $R$-module spectrum $X^\zeta$ should be thought of as a twisted form of $X\wedge R$.
The above construction is functorial in the following sense.
Given a map $R\to S$ of $A_\infty$-ring spectra,
one gets a corresponding map
\[
f:B\GL_1(R)\to B\GL_1(S).
\]
If $\zeta: X\to B\GL_1(R)$ is as above, then one finds that
\[
X^{f\circ\zeta}=X^\zeta\wedge_R S.
\]

Now consider the $J$-homomorphism $O\to \GL_1(\SS)$ that
sends a linear automorphism of a vector space to the corresponding self-homotopy equivalence
of its one-point compactification.
Its delooping
\[
J:BO\to B\GL_1(S)
\]
associates to a vector bundle $V$ a spherical fibration $\mathit{Sph}(V)$.
Let $\iota:B\GL_1(\SS)\to B\GL_1(R)$ be induced by the unit map $\SS\to R$.
If $V$ is a vector bundle and $\zeta:=\mathit{Sph}(V)$ the corresponding spherical fibration,
then nullhomotopies of $\iota\circ \zeta$ correspond to $R$-orientations of $V$:
\[
\xymatrix{
BO\ar[r]^(.4)J&B\GL_1(\SS)\ar[r]^\iota &B\GL_1(E)\\
X\ar[u]^V\ar[ur]_\zeta
}
\]
Indeed, $X^\zeta$ is the usual Thom spectrum of $V$, and
$X^{\iota\circ \zeta}=X^\zeta\wedge R$ is the spectrum that we want to trivialize.
A homotopy $\iota\circ \zeta\simeq *$ induces an isomorphism $X^\zeta\wedge R\simeq X\wedge R$.

Now suppose that we want to functorially $R$-orient every vector bundle. 
Then we would need to chose a nullhomotopy for the composite
\[
BO\xrightarrow{J} B\GL_1(\SS)\xrightarrow{\iota} B\GL_1(R).
\]
If we only wanted to functorially $R$-orient spin bundles, then we would need a nullhomotopy of
\begin{equation}\label{bssr}
B\mathit{Spin}\to BO\xrightarrow{J} B\GL_1(\SS)\xrightarrow{\iota} B\GL_1(R).
\end{equation}
Similarly, if we wanted to $R$-orient string bundles, then we would need a nullhomotopy of
\[
BO\langle8\rangle\to BO\xrightarrow{J} B\GL_1(\SS)\xrightarrow{\iota} B\GL_1(R).
\]
At this point, finding those nullhomotopies is still a rather hard problem.
For example, if $R$ was $KO$-theory, and if we were trying to construct the Atiyah-Bott-Shapiro orientation (\ref{MKO}),
then we wouldn't be able to handle that.
The reason is that the space
\[
\Map \big(B\mathit{Spin},B\GL_1(KO)\big)
\]
is very big. 
It is hard to tell which map $B\mathit{Spin}\to B\GL_1(KO)$ one is looking at.
And in particular, it is hard to tell if a map is null.

To solve that problem, we impose more conditions: we ask that (\ref{bssr}) be nullhomotopic through $E_\infty$ maps!
Of course, that condition doesn't make any sense yet, because we had only assumed that $R$ was $A_\infty$.
But if $R$ is an $E_\infty$-ring spectrum, then $\GL_1(R)$ is an $E_\infty$-space
and the condition does make sense.

So let's assume that $R$ is an $E_\infty$, and let $\gl_1(R)$ be the spectrum associated to $\GL_1(R)$.
Since $\gl_1(R)$ is obtained by delooping an $E_\infty$-space, it is necessarily $(-1)$-connected
(there are also ways of adding negative homotopy groups to $\gl_1(R)$, but that's not relevant for the present discussion).

Let $Y$ be a $(-1)$-connected spectrum and let $X:=\Omega^\infty Y$ be its zeroth space.
We can then consider infinite loop maps
\[
\zeta:X\to B\GL_1(R),
\]
or equivalently, maps of spectra
\[
\zeta^\infty: Y\to \Sigma\,\gl_1(R).
\]
Since $\zeta$ is an infinite loop map, the Thom spectrum $X^\zeta$ is an $E_\infty$-ring spectrum.
$E_\infty$-orientations then correspond to nullhomotopies of $\zeta$ through infinite loop maps.
Equivalently, they correspond to nullhomotopies of $\zeta^\infty$.
In our case of interest, we see that
$E_\infty$ maps $M\mathit{Spin}\to KO$ correspond to nullhomotopies of the composite
\begin{equation}\label{bbSS}
b\mathit{spin}\,\to\, bo\,\xrightarrow{J}\,\Sigma\,\gl_1(\SS)\,\to\, \Sigma\,\gl_1(KO).
\end{equation}
Similarly, $E_\infty$ maps $M\mathit{String}\to \tmf$ correspond to nullhomotopies of
\[
bo\langle8\rangle \,\to\, bo\,\xrightarrow{J}\,\Sigma\,\gl_1(\SS)\,\to\, \Sigma\,\gl_1(\tmf).
\]

A nullhomotopy of a composite $f\circ g$ is the same thing as an extension of $f$ over the cone of $g$.
So a nullhomotopy of (\ref{bbSS}) is equivalent to a dotted arrow
\begin{equation}\label{sbgg}
\begin{matrix}
\xymatrix{
\Sigma^{-1}b\mathit{spin}\ar[r]^(.52)J&\gl_1(\SS)\ar[r]\ar[d]&C\ar@{.>}[dl]\\
&\gl_1(KO)
}
\end{matrix}
\end{equation}
making the diagram commute. Similarly, dotted arrows
\begin{equation}\label{sbg2}
\begin{matrix}
\xymatrix{
\Sigma^{-1}bo\langle8\rangle\ar[r]^(.52)J&\gl_1(\SS)\ar[r]\ar[d]&D\ar@{.>}[dl]\\
&\gl_1(\tmf)
}
\end{matrix}
\end{equation}
correspond to $E_\infty$ orientations $M\mathit{String}\to \tmf$.
The horizontal lines in (\ref{sbgg}) and (\ref{sbg2}) are cofiber sequences, and we have desuspended all our spectra to simplify the notation.
The set of extensions (\ref{sbg2}) is either empty, or a torsor over the group $\big[bo\langle8\rangle,\gl_1(\tmf)\big]$.

From now on, we pick a prime $p$, and assume implicitly that all spectra are $p$-completed.
Fix $n\ge 1$. Bousfield and Kuhn constructed a functor $\Phi$ from spaces to spectra, that factors $K(n)$-localization as
\[
\xymatrix{
\text{Spectra}\ar[rr]^{L_{K(n)}}\ar[dr]_{\Omega^\infty}&&\text{Spectra}.\\
&\text{Spaces}\ar[ur]_\Phi
}
\]
Apart from the difference of $\pi_0$, 
the zeroth space $\Omega^\infty(\gl_1(R))=GL_1(R)$ of the spectrum $\gl_1(R)$ is the same as the zeroth space of $R$.
Since $L_{K(n)}$ kills Eilenberg-McLane spectra, the difference of $\pi_0$ doesn't matter, and so we get
\[
L_{K(n)}\big(\gl_1(R)\big)\,\simeq\, L_{K(n)}(R).
\]
More generally, if $X$ and $Y$ are spectra such that for some $m\ge 0$, 
the $m$-th connected cover of $\Omega^\infty X$ agrees with the the $m$-th
connected cover of $\Omega^\infty Y$, then $L_{K(n)}(X)\,\simeq\, L_{K(n)}(Y)$.

Since $KO$ is $K(1)$-local, the localization map $\gl_1(KO)\to L_{K(1)}\big(\gl_1(KO)\big)$
induces a map
\[
L:\gl_1(KO)\to KO.
\]
The spectrum $\gl_1(KO)$ being connected, there is no hope for $L$ to be an isomorphism.
But Bousfield proved that it is a $\pi_*$-isomorphism for $*>2$.
Going back to (\ref{sbgg}), we note that the first non-zero homotopy group of $\Sigma^{-1}b\mathit{spin}$ is in dimension $3$. 
That is exactly the range above which $\gl_1(KO)$ looks like $KO$.
So the obstruction to constructing our orientation (\ref{MKO}) may be taken to be the composite
\[
\Sigma^{-1}b\mathit{spin}\to \gl_1(KO)\xrightarrow{L} KO.
\]
It lives in $[\Sigma^{-1}b\mathit{spin},KO]=KO^1(b\mathit{spin})$, which can be shown to be zero.
The calculation
\[
KO^1\big(b\mathit{spin}\big)=0
\]
is actually not too hard: 
it follows from the knowledge of the cohomology operations in $KO$-theory, which is something that one computes using Landweber exactness.
So one gets the existence of a dotted arrow in (\ref{sbgg}).

The above discussion was done at the cost of completing at a prime.
To do things correctly, one should complete at each prime individually, do something rationally, 
and then assemble the results using Sullivan's arithmetic square.
Eventually, one sees that the map $\Sigma^{-1}b\mathit{spin}\to \gl_1(KO)$ is null (no completion any more),
With a little bit more work, one can parametrize the set of nullhomotopies of (\ref{bbSS}),
which amounts to parametrizing the $E_\infty$-orientations (\ref{MKO}). 

We now proceed to the case of $\tmf$.
That's a more complicated story, but it also leads to something more interesting.
First of all, one can generalize Bousfield's result and show that the fiber of the sequence 
\[
F\to \gl_1(\tmf)\xrightarrow{L} L_{K(1)\vee K(2)}\big(\gl_1(\tmf)\big)
\]
has no homotopy groups in dimensions $3$ and above [I corrected the bound stated in the talk].
So, as far as mapping $\Sigma^{-1}bo\langle8\rangle$ is concerned, we may replace $\gl_1(\tmf)$ by its
$K(1)\vee K(2)$-localization.
As a consequence, if we want to produce a nullhomotopy for our map
\[
\Sigma^{-1}bo\langle8\rangle\to \gl_1(\tmf),
\]
we may as well produce one for the composite
\[
\Sigma^{-1}bo\langle8\rangle\to \gl_1(\tmf)\xrightarrow{L} L_{K(1)\vee K(2)}\big(\gl_1(\tmf)\big).
\]
 
Since we understand better the localizations at individual Morava $K$-theories, we consider the Hasse pullback square
\begin{equation}\label{g4l}
\begin{matrix}
\xymatrix{
\gl_1(\tmf)\ar@{-}[r]^(.35){\substack{\text{in}\\ \text{dim}\,\ge\,\, 4\\{\phantom{\big |}|\phantom{\big |}}\\ \textstyle\simeq}}&L_{K(1)\vee K(2)}\big(\gl_1(\tmf)\big)\ar[r]\ar[d]&L_{K(2)}\big(\gl_1(\tmf)\big)\ar[d]\\
&L_{K(1)}\big(\gl_1(\tmf)\big)\ar[r]&L_{K(1)}L_{K(2)}\big(\gl_1(\tmf)\big).
}
\end{matrix}
\end{equation}
The spectrum $\Sigma^{-1}bo\langle8\rangle$ doesn't have any $K(2)$ cohomology, 
and therefore
\[
\Big[\Sigma^{-1}bo\langle8\rangle,L_{K(2)}\big(\gl_1(\tmf)\big)\Big]_*=0.
\]
So as far as mapping $\Sigma^{-1}bo\langle8\rangle$ into it,
the square (\ref{g4l}) behaves like a fiber sequence, and we get a long exact sequence
\[
\Big[\Sigma^{-1}bo\langle8\rangle,\gl_1(\tmf)\Big]\to
\Big[\Sigma^{-1}bo\langle8\rangle,L_{K(1)}\big(\gl_1(\tmf)\big)\Big]\to
\Big[\Sigma^{-1}bo\langle8\rangle,L_{K(1)}L_{K(2)}\big(\gl_1(\tmf)\big)\Big]\cdots
\]
Now, we wish to apply Bousfield and Kuhn's result to identify the above spectra.
At first glance, things look pretty good because both $L_{K(1)}(\gl_1(\tmf))$
and $L_{K(1)}L_{K(2)}(\gl_1(\tmf))$ are $K(1)$-local, and hence ``look a lot like $K$-theory''.
We also saw that in that range of dimensions [what range of dimensions?],
there are no maps from $\Sigma^{-1}bo\langle8\rangle$ into $L_{K(1)}(\gl_1(\tmf))$ and into $L_{K(1)}L_{K(2)}(\gl_1(\tmf))$...
But there are some subtelties.

If we apply the argument with the Bousfield Kuhn functor to the diagram (\ref{g4l}), we get a square 
\begin{equation}\label{gKKK}
\begin{matrix}
\xymatrix{
\gl_1(\tmf)\ar[r]\ar[d]&L_{K(2)}(\tmf)\ar[d]\\
L_{K(1)}(\tmf)\ar[r]&L_{K(1)}L_{K(2)}(\tmf)
}
\end{matrix}
\end{equation}
that is a pullback square above dimension 4.
The above diagram is very similar to the Hasse square for $\tmf$.
The only difference is that it has $\gl_1(\tmf)$ insted of $\tmf$ in its upper left corner.
The right vertical map
\[
L_{K(2)}(\tmf)\to L_{K(1)}L_{K(2)}(\tmf)
\]
is simply the localization map $X\to L_{K(1)}(X)$ applied to $X=L_{K(2)}(\tmf)$,
but the lower map
\begin{equation}\label{Atkin}
L_{K(1)}(\tmf)\to L_{K(1)}L_{K(2)}(\tmf)
\end{equation}
is much more subtle.
In some sense, the whole game of producing the string orientation, is to understand that map.

At the beginning, we had a few of ad-hoc ways of understanding that map.
Later, Charles Rezk found extremely beautiful formula for it. 
It says that, at the level of homotopy groups, the map (\ref{Atkin}) is given by $1-U_p$, where $U_p$ is the Atkin operator.
Knowing that fact, we may replace (\ref{gKKK}) by
\begin{equation}\label{gKKK2}
\begin{matrix}
\xymatrix@C=1.5cm{
\gl_1(\tmf)\ar[r]\ar[d]&\tmf\ar[d]\\
L_{K(1)}(\tmf)\ar[r]^{\,\,1-U_p\,\,}&L_{K(1)}(\tmf).
}
\end{matrix}
\end{equation}
The above square ends up giving us enough understanding of the homotopy type of $\gl_1(\tmf)$
[we refer the reader to \cite{AHR07} for the actual computations].

\[
\sim\qquad\sim\qquad\sim
\]\phantom{.}

Let me end this talk by emphasizing the number of really amazing things that come out of the square (\ref{gKKK2}).
First of all, it gives you the string orientation.
The question of orientation boils down to the question of understanding the homotopy type of $\gl_1(\tmf)$;
using the square (\ref{gKKK2}) and the description of its bottom arrow, one can then make the required calculations.
Secondly, the homotopy fiber of $1-U_p:L_{K(1)}(\tmf)\to L_{K(1)}(\tmf)$ is an extremely interesting spectrum.
One can describe its homotopy groups in terms of modular forms,
but if one tries to compute them explicitly, one encounters some unsolved problems in number theory.
Finally, we conjecture that there is a relationship between that homotopy fiber and 
[smooth structures on free loop spaces of spheres].

\bibliographystyle{plain}


\end{document}